\documentclass{article}
\usepackage{amsfonts, amssymb, latexsym, graphicx}
\usepackage[all]{xy}
\newtheorem{theorem}{Theorem}
\newtheorem{lemma}{Lemma}

\begin{document}
\title{Totally real pencils of cubics with respect to sextics}
\author{S\'everine Fiedler-Le Touz\'e}
\maketitle
\begin{abstract}
A real algebraic plane curve $A$ is said to be dividing if its real part $\mathbb{R}A$ disconnects its complex part $\mathbb{C}A$. A pencil of curves is totally real with respect to $A$ if it has only real intersections with $\mathbb{C}A$. If there exists such a pencil, then $A$ is dividing, this is the case for the M-curves. Can conversely any dividing curve be endowed with a totally real pencil? We study here the case of M-2-sextics having $2$ or $6$ empty exterior ovals and one non-empty oval surrounding $6$ or $2$ empty ovals. Such sextics are always dividing. We prove that they may actually be endowed with a totally real pencil of cubics.  
\end{abstract}

\section{Introduction}
The aim of this note is to prove the following theorem:
\begin{theorem} 
i) Any M-2-sextic with real scheme $\langle 2 \amalg 1 \langle 6 \rangle \rangle$ or $\langle 6 \amalg 1 \langle 2 \rangle \rangle$ may be endowed with a totally real pencil of cubics. It suffices to choose eight base points distributed on the eight empty ovals.

ii) The empty ovals of an M-2-sextic with real scheme $\langle 2 \amalg 1 \langle 6 \rangle \rangle$ lie in convex position.  Let $1, \dots 8$ be eight generic points distributed in successive empty ovals, such that $8$ and $1$ lie in the two outer ovals and the oval containing $7$ forms a negative pair with the non-empty oval $O$.  Then, the successive distinguished cubics of the pencil determined by $1, \dots 8$ are as shown in Figure~1, where $\{ X, Y \} = \{ 2, 9 \}$ ($9$ is the ninth base point), and $O$ is represented with a dotted line. If $2$ is chosen on the empty oval instead of inside, then the pencil is totally real. 

iii) Along a rigid isotopy between two sextics  $\langle 2 \amalg 1 \langle 6 \rangle \rangle$, this totally real pencil of cubics may be deformed continuously, preserving total reality, with at most one type of degeneration, after which the positions of the two base points $2$ and $9$ on the distinguished cubics are swapped.
\end{theorem}

\begin{figure}
\centering
\includegraphics{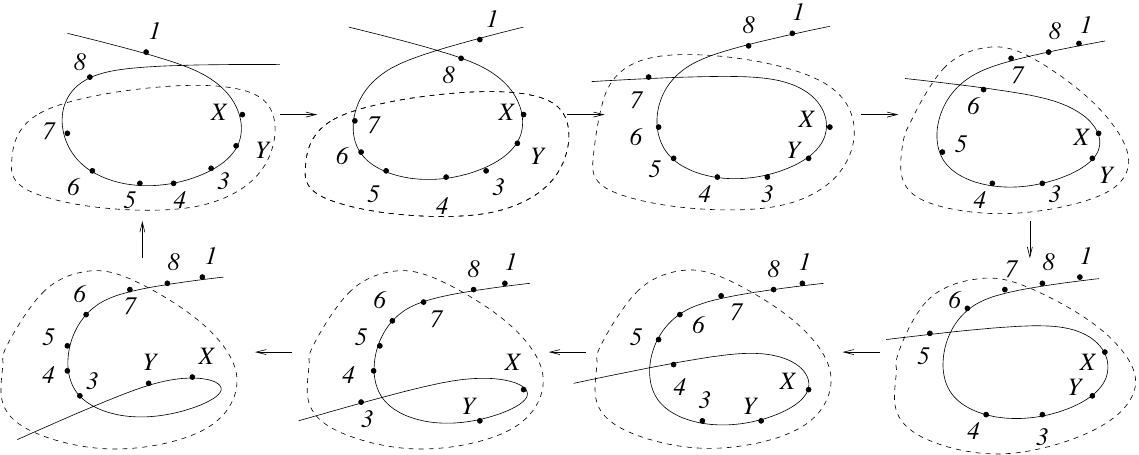}
\caption{\label{totalrealpen} Pencil of cubics with eight base points distributed in the eight empty ovals of $\langle 2 \amalg 1 \langle 6 \rangle \rangle$} 
\end{figure}   

Let us recall some definitions and some facts before stepping forward to the proof.
Let $A$ be a real algebraic non-singular plane curve of degree $m$. Its complex part $\mathbb{C}A$ is a Riemannian surface of genus $g = (m-1)(m-2)/2$; its real part $\mathbb{R}A$ is a collection of circles embedded in $\mathbb{R}P^2$. The maximal number of circles is $M = g+1$, a curve realizing this upper bound is called
M-curve. A circle embedded in $\mathbb{R}P^2$ is called oval 
or pseudo-line depending on whether it realizes the class $0$ or $1$ of $H_1(\mathbb{R}P^2)$. If $m$ is even, the curve has only ovals, otherwise it has
ovals plus one unique pseudo-line, also called odd component. The real scheme of $A$ is the isotopy type of $\mathbb{R}A \subset \mathbb{R}P^2$, it may be encoded with Viro's notation that we explain here using a simple example: take a sextic with $\alpha$ exterior ovals plus a non-empty oval surrounding $\beta$ empty ovals. Its real scheme will be denoted by $\langle \alpha \amalg 1 \langle \beta \rangle \rangle$. 
For a curve of odd degree, the odd component will be encoded by $\mathcal{J}$.
The complex conjugation acts on  $\mathbb{C}A$ with $\mathbb{R}A$ as fixed points set. Thus $\mathbb{C}A \setminus \mathbb{R}A$ is connected or splits into two homeomorphic halve that are exchanged by $conj$, in this case we say that $A$ is dividing, or of type I. Otherwise, $A$ is of type II. For a dividing curve, the orientation of any half of $\mathbb{C}A \setminus \mathbb{R}A$ induces an orientation on the real components, called complex orientation. The complex orientation, figured by arrows on the real circles, is defined only up to complete reversion. 
If the orientations of two nested ovals induce an orientation of the annulus they bound, then we say that these two ovals form a positive pair; otherwise they form a negative pair. The numbers of positive and negative pairs are denoted by $\Pi_+$
and $\Pi_-$. If $A$ has odd degree, then each oval $O$ may be provided with a sign: consider the M\"obius band $\mathcal{M}$ obtained by cutting away the interior of $O$ from $\mathbb{R}P^2$. The classes $[O]$ and $[2\mathcal{J}]$ of $H_1(\mathcal{M})$ coincide or are opposite. In the first case, we say that $O$ is negative; otherwise $O$ is positive.
Let $\Lambda_+$ and $\Lambda_-$ be the numbers of positive and negative ovals.
A real scheme of given degree is of type I if it is realizable by dividing curves only; it is of type II if it is realizable by non-dividing curves only. Otherwise, it is indefinite. 

The M-curves are clearly dividing. The so-called hyperbolic curves are also dividing. A hyperbolic curve of degree $m = 2k$ or $2k+1$, consists in $k$ nested ovals, plus one pseudo-line if $m$ is odd. A pencil of lines whose base point is chosen in the innermost oval sweeps out the curve in such a way that the $m$ intersections are always real.
One says that this pencil of lines is totally real with respect to the hyperbolic curve. Starting from this observation. Rokhlin \cite{R} presents a beautiful argument proving that if an algebraic curve is swept out by a totally real pencil of lines, then this curve is dividing. The argument generalizes to pencils of curves of higher degrees. Can conversely any dividing curve be endowed with some totally real pencil?
A weaker conjecture suggested implicitely in \cite{R} is that any curve whose real
scheme is of type I may be endowed with a totally real pencil.
It turns out that the M-curves may indeed be endowed with suitable pencils of degree $m-2$, see \cite{G}, page 348. 
We shall consider always generic real pencils with only real base points,  the intersections are counted with multiplicities. A homogeneous polynomial of degree $d$ in three variables has $d(d+3)/2+1$ monomials, hence the space of real algebraic curves of degree $d$ is a projective space $\mathbb{R}P^N$, with $N = d(d+3)/2$.
The discriminantal hypersurface $\Delta$ formed by the singular curves has degree $3(d-1)^2$. A pencil of curves of degree $d$ has $d^2$ base points, any $N-1$ of them determine the pencil, and the number of singular curves is $3(d-1)^2$.  
Two non-singular curves of degree $d$ are rigidly isotopic if they are in the same
connected component of $\mathbb{R}P^N \setminus \Delta$.

Let us give a qualitative description of pencils with low degree.
A pencil of conics has four base points and three singular conics (double lines).
A pencil of cubics \cite{FLT2} has nine base points, and is determined by any eight of them. Exactly eight of the twelve singular cubics are {\em distinguished cubics\/} that is to say real nodal cubics with some base points on the loop. The circle parametrizing the real part of the pencil is divided in eight portions by the distinguished cubics, four pairwise non-consecutive portions contain only smooth cubics with an oval passing through an even number of base points, the other four portions contain connected cubics, plus possibly subportions bounded by a pair of singular cubics, one with an isolated node, the other with a loop passing through no base point. The inside of the subportion consist in cubics with an oval passing through no base point.  
Let us present some examples of totally real pencils. For an M-quartic $\langle 4 \rangle$, a pencil of conics, with the four base points distributed in the four ovals
will do. Similarly, the curve of degree $8$ with real scheme 
$\langle 1 \langle 1 \rangle \amalg 1 \langle 1 \rangle \amalg 1 \langle 1
\rangle \amalg 1 \langle 1 \rangle \rangle$ formed by four nests of depth $2$ lying outside one another admits a totally real pencil of conics, it suffices to choose the base points inside of the innermost ovals. For an M-quintic $\langle \mathcal{J} \amalg 6  \rangle$, one finds a suitable pencil of cubics with six base points distributed on the six ovals, and two further chosen on the odd component $\mathcal{J}$. As this component must cut any cubic an odd number of times, the required $15$ real intersections are granted.  
 
By a congruence due to Kharlamov \cite{K-V}, the real schemes $\langle 2 \amalg 1 \langle 6 \rangle \rangle$ or $\langle 6 \amalg 1 \langle 2 \rangle \rangle$ are both of type I. We use this fact in the proof, so we confirm that the sextics with these two schemes do not contradict Rokhlin's conjecture. Rokhlin claimed that he could prove the very same statement {\em without\/} using the fact that the sextics are dividing. It's a stronger result. Unfortunately, his proof was never published and is now lost. Nikulin and Kharlamov \cite{N}, \cite{K} proved that the rigid isotopy class of a non-singular sextic is determined by its real scheme, plus its type I or II. Thus, any two M-2-sextics realizing the same scheme among our two, are rigidly isotopic.  

The fact that the empty ovals of a sextic $\langle 2 \amalg 1 \langle 6 \rangle \rangle$ lie in convex position allows us to find a precise description of the pencil, using the tools developed in \cite{FLT2}. In this paper, we have classified pencils of cubics with eight base points lying in convex position.
Let us recall briefly some ideas developed there. Consider eight points $1, \dots 8$ lying in convex position in the plane, one may associate to these points
two combinatorial objects called respectively {\em list\/} and {\em combinatorial pencil\/}. The combinatorial pencil $\mathcal{P}(1, \dots 8)$ is the list of topological types (cubic, base points) realized by the eight successive distinguished cubics
of the actual pencil.  
The list $L(1, \dots, 8)$ is the list of the $56$ conics passing through $5$ of the points, enhanced for each conic, with the position of each of the remaining $3$ point (inside or outside). How many lists may be realized when one lets the points move? Written in full extent, the lists are full of redundancies. Using rational cubics, we found an encoding avoiding redundancies, and managed to classify the lists. 
Up to the action of $D_8$, the number of lists realizable by eight points lying in stricly convex position, no six on a conic, is $49$. 
We considered first configurations of seven points $1, \dots 7$ in convex position and the list of seven rational cubics passing through them, with node at one of them. For convenience, we call shortly cubic a topological type (cubic, seven points). There are $14$ possible lists of seven cubics, we denote these lists by $n\pm$, $n = 1, \dots 7$, they are all equivalent up to cyclic permutation and symmetries, see Proposition 1 and Figure~3 in \cite{FLT2}. The cubics of $n\pm$ with node at $n$ have a loop passing through all other six points. The encoding may be generalized to configurations of eight points, one has eight lists $\hat n, n = 1, \dots 8$ of seven rational cubics passing through all of the points but $n$. It turns out that the data $L(1, \dots 8)$ and $(\hat 1, \dots \hat 8)$ are equivalent. A list is {\em maximal\/} if one point lies outside of all the conics determined by five others. 
Move the eight points preserving the strict convexity. An {\em elementary change\/} is the change induced on a list by a move letting one point cross a conic through five others. Up to the action of the dihedral group $D_8$, there are $19$ elementary changes. Four of the $49$ lists are {\em nodal\/}, that is to say realizable by eight points on a nodal cubic, one of them being the node. 
Now, let us watch the connections between lists and pencils. If one moves the eight points preserving the strict convexity, the first degeneration of the (actual) pencil occurs when one of the points $1, \dots 8$ comes together with $9$, or when six of the points $1, \dots 8$ become coconic, $9$ is then aligned with the other two. At this moment, two singular cubics of the pencil come together to yield one singular cubic with node at one of the points $1, \dots 8$, or a reducible cubic.  
Bezout's theorem allows to recover the initial  pair of (combinatorial) singular cubics. 
The combinatorial distingushed cubics are encoded as explained in page 21 of \cite{FLT2}. Among the $19$ elementary changes, only one is {\em inessential\/} that is to say, it doesn't alter the combinatorial pencil. If a pencil fits to two different
lists, then these lists differ by a non-essential change. Any non-nodal list corresponds to one single combinatorial pencil, a nodal list gives rise to several pencils that are obtained successively from one another by swaps of $9$ with other base points. Up to the action of $D_8$, the number of combinatorial pencils realizable with eight base points in convex generic position is $45$.

\section{Proofs}
\begin{lemma}
Consider a real sextic $C_6$ with one non-empty oval $O$.
There is a natural cyclic ordering of the empty ovals of $C_6$, consistant with the ordering given by the complete pencils of lines based at the inner ovals and the partial pencils of lines based at outer ovals and sweeping out $O$. For this cyclic ordering, the set of empty ovals splits in two consecutive chains (inner, outer). Moreover, the inner ovals lie in convex position inside of $O$. If $C_6$ has at most
two outer ovals, then all of the empty ovals lie in convex position.
\end{lemma}

{\em Proof of the Lemma:\/}
Given an empty oval $X$ of $C_6$ we will often have to consider one point chosen in the interior of $X$. For simplicity, we call this point also $X$, it will be clear from the context whether we speak of the oval or of the point. 
Let $A, B$ be two outer ovals, $C$, $D$ be two inner ovals. By Bezout's theorem with the line $CD$, one segment $[CD]$ is entirely contained in $O$. Assume that the line $AB$ cuts this segment, and consider the pencil of conics based at $A, B, C, D$. The conics of this pencil cut the set of ovals $A, B, C, D, O$ at $12$ points, so the other empty ovals cannot be swept out, this is a contradiction. So the line $AB$ must cut the other segment of line $[CD]'$. Let now $C_6$ have more than two inner ovals, and let $E$ be a third inner oval. The points $C, D, E$ determine four triangles in $\mathbb{R}P^2$, one of them is entirely contained in $O$, let us call it {\em principal triangle\/} $CDE$. Let $F$ be another empty oval. Using a pencil of conics based at $C, D, E, F$, we prove that $F$ cannot lie inside of the principal triangle $CDE$.  Assume that two outer ovals $A$ and $B$ lie in two different (non-principal) triangles $CDE$. Then the conic through $A, B, C, D, E$ cuts the set of ovals $A, B, C, D, E, O$ at $14$ points, this is a contradiction. $\Box$

Note also that the cyclic ordering gives rise to a choice of a {\em principal segment\/} $[XY]$ connecting two consecutive ovals $X, Y$, if $X, Y$ are inner ovals, this segment lies inside of $O$. 

{\em Proof of i)\/}
Let us consider the M-2-sextics $\langle 2 \amalg 1 \langle 6 \rangle \rangle$
and $\langle 6 \amalg 1 \langle 2 \rangle \rangle$, and a pencil of cubics
with eight base points distributed on the eight empty ovals, all of the cubics of this pencil cut the set of empty ovals at $16$ points. Assume that the pencil
is not totally real, then it has a cubic $C_3$ that does not cut $O$. This cubic is
maximal, its oval $\mathcal{O}$ passes through the inner ovals of $C_6$ and its odd component $\mathcal{J}$ passes through the outer ovals of $C_6$ as shown in Figure~2. 

\begin{figure}
\centering
\includegraphics{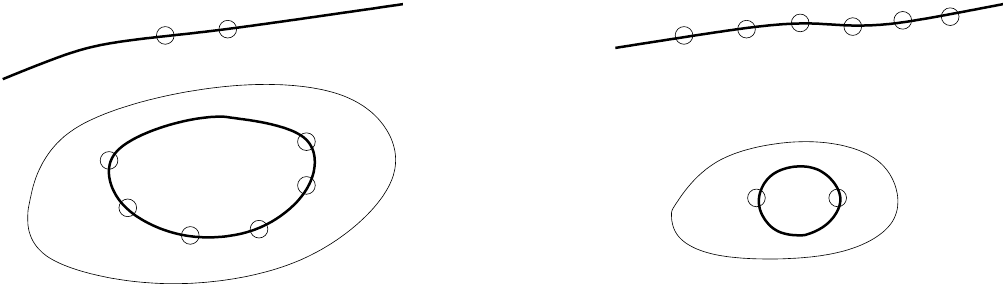}
\caption{\label{problemcub} The cubic $C_3$} 
\end{figure}   

Note that the ninth base point lies on the odd component, outside from $O$. The cyclic orderings of the inner ovals on $\mathcal{O}$ and of the outer ovals on $\mathcal{J}$ are consistant with the natural cyclic ordering defined hereabove, this
follows from Bezout's theorem between $C_3$ and the pencils of lines based at inner ovals. In an affine plane where the line at infinity cuts none of the principal
segments $[XY]$ (as in Figure~2), the naturally ordered empty ovals have alternating orientations by Fiedler's theorem \cite{F2}. Let us move in the pencil of cubics, starting from $C_3 = C_3(0)$ in one of the two possible directions. In the upper part of Figure 3, the cubic $C_3$ is drawn in bold, and the dotted oval is the oval $\mathcal{O}(t)$ of some cubics $C_3(t)$ close to $C_3$. During the motion, two mobile arcs with fixed endpoints, one arc $s$ on the oval $\mathcal{O}(t)$,  the other $s'$ on the odd component $\mathcal{J}(t)$, move towards each other until they are glued with a non-isolated double point. Let us denote by $C_3(1)$ the corresponding singular cubic. For the sextic with six inner ovals, let us denote the base points of the pencil by $1, \dots 8$ as indicated in Figure~3. The mobile arcs $27$, $56$ and $34$ of $\mathcal{O}(t)$ are outside from the oval $\mathcal{O}$ of the initial cubic $C_3$.  Thus $s$ is one of these three arcs. Note that each mobile arc of $\mathcal{O}(t)$ is homotopic with fixed extremities to a segment of line contained in $O$. The line $68$ cuts the segment of lines $[27]$ interior to $O$, so this line cuts $C_3(t)$ at $6, 8$ and a point of the arc $27$. During the motion, the mobile arc $56$ must stay inside of a zone delimited by the lines $68$, $54$, and the arc $56$ of the initial oval $\mathcal{O}$. As $\mathcal{J}(t)$ may a priori intersect only the side $54$, and this only once, this zone is not cut by $\mathcal{J}(t)$, see the upper part of Figure~3. Symmetrically, the mobile arc $34$ is bound to stay inside of a zone that is not cut by $\mathcal{J}(t)$ either. Therefore, the singular cubic $C_3(1)$ is obtained attaching the arc $27$ of $\mathcal{O}(t)$ to some arc of the odd component. 
The singular cubic $C_3(1)$ cuts $C_6$ at $18$ points. If one provides both curves with a complex orientation, and perturb all double points of their product according to these orientations, one gets a curve $C_9$ of degree $9$ with $l = 10$ ovals. This curve is dividing, with a complex orientation induced by that of the initial curves, see \cite{F2}. Rokhlin-Mishachev's formula \cite{R} yields: $2(\Pi_+ - \Pi_- ) + \Lambda_+ - \Lambda_- = l - k(k+1) = -10$. The $n = 2$ or $6$ base points chosen on the outer ovals divide the odd component $\mathcal{J}(t)$ in $n$ arcs, one of them contains $s'$ or is the whole of $s'$. Clearly, the data $\Pi_\pm$, and $\Lambda_\pm$ depend neither on the actual position of $s'$ nor on the position of the double point of $C_3(1)$ inside or ouside of $O$. 
Up to a symmetry with respect to a vertical axis passing through the middle of the
Figure, changing $C_3(1)$, the orientation of $O$ may be chosen at leisure.
In the middle part of Figure~3, we have drawn possible cubics $C_3(1)$ for either sextic with the non-empty oval $O$ enhanced with some orientation. In the bottom part, we have represented the two curves $C_9$ obtained.  One gets  $2(\Pi_+ - \Pi_-) + \Lambda_+ - \Lambda_- = -8$ in both cases, this is a contradiction.  
$\Box$ 

\begin{figure}
\centering
\includegraphics{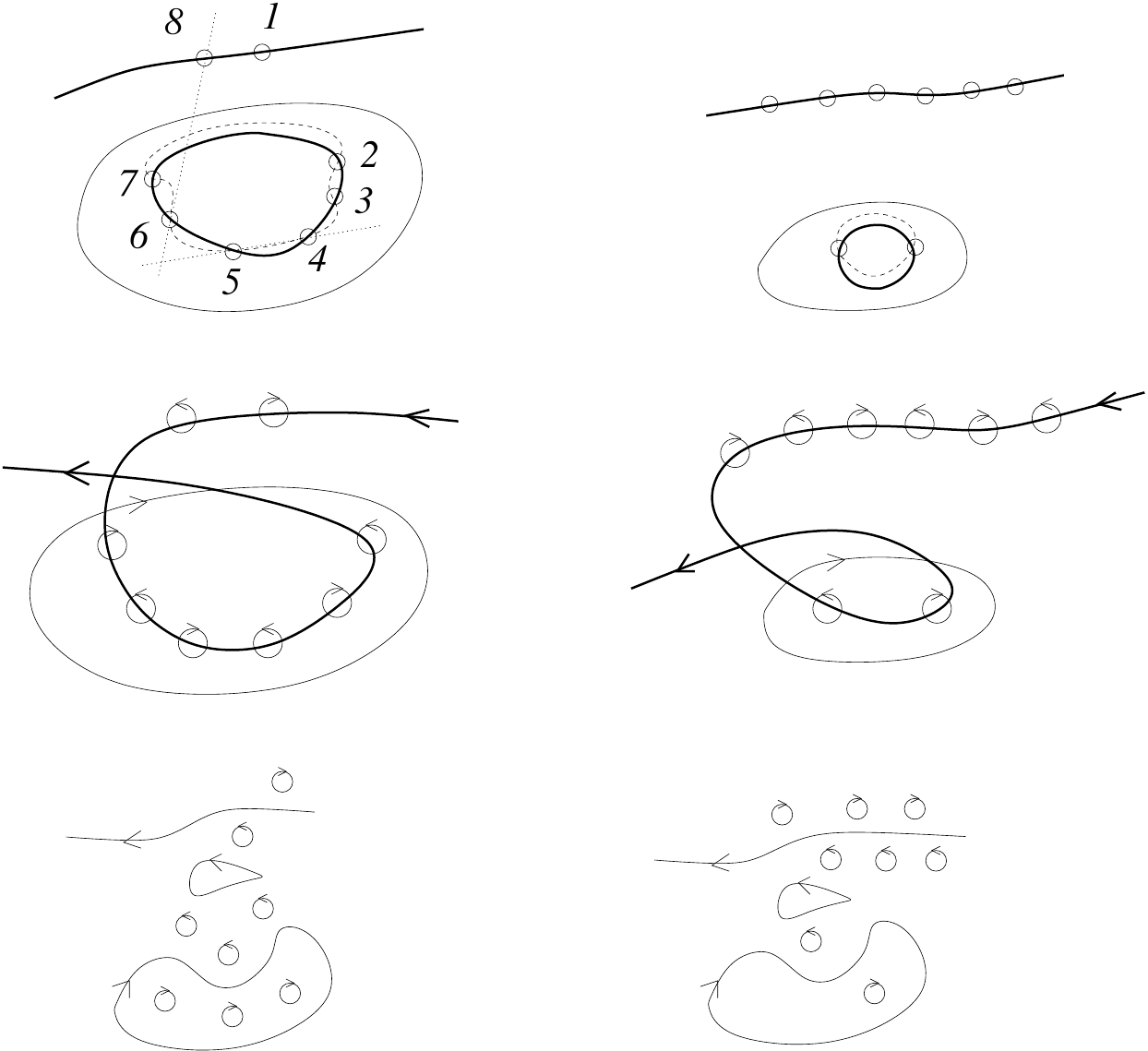}
\caption{\label{movetosing} The next singular cubic in the pencil leads to a contradiction} 
\end{figure}   

{\em Proof of ii) and iii)\/} For this part, we use the notations from \cite{FLT2}. 
Let us determine the list $L(1, \dots 8)$ realized by eight generic points distributed in the eight empty ovals of the sextic $\langle 2 \amalg 1 \langle 6 \rangle \rangle$ as indicated in the upper part of Figure~3. By Bezout's theorem, one has 
$\hat 8 = 1\pm$, $\hat 1 = 8\pm$, and $\hat n \in \{ 8\pm, 1\pm \}$
for $n = 2, \dots 7$.
On the other hand, if $\hat 8 = 1+$, then $\hat 1 \in \{ 8+, 6\pm, 4\pm, 2\pm \}$, see Figure~39 in \cite{FLT2}.
So the list $L(1, \dots 8)$ has either a pair $\hat 8 = 1+$, $\hat 1 = 8+$, or a
pair $\hat 8 = 1-$, $\hat 1 = 8-$ (the second pair is deduced from the first by the action of the element $(+1)(48)$ of $D_8$). Figure~9 from \cite{FLT2} shows the $32$ lists with $\hat 8 = 1+$, $\hat 1 = 8+$. Only two of them satisfy the
condition on the $\hat n$, $n = 2, \dots 7$, namely the list $L_1 = \max(\hat 8 = 1+)$ and $L_{32} = \max(\hat 1 = 8+)$. The notation $\max(\hat 8 = 1+)$ designs the maximal list such that $8$ lies outside of all the conics determined by other five points, and $\hat 8 = 1+$. Up to the action of $(+1)(48)$, let us assume
that $L(1, \dots 8)$ is $\max(\hat 1 = 8+)$ or $\max(\hat 1 = 8-)$.
These two lists are nodal, $\max(\hat 1 = 8-)$ gives rise to nine combinatorial pencils (Figures~33 and 34 in \cite{FLT2}) and $\max(\hat 1 = 8+)$ is realizable with the first three combinatorial pencils among those nine. (The two lists differ by a non-essential elementary change.) Let us choose a cubic $C_3$ of the pencil, provide it with an orientation, perturb the union $C_3 \cup C_6$ and check Rokhlin-Mishachev's formula for the obtained curve $C_9$.  Assume that one of the first seven pencils of Figure~34 is realized,  choose $C_3$ to be $(1-, 3)$. We find a contradiction for both orientations of the non-empty oval $O$ of $C_6$. Thus, the pencil is one of the last two. Let us now choose $C_3$ to be the cubic $(1+, 12)$, the case where $O$ and the oval surrounding $7$ form a positive pair yields again a contradiction. One checks that the other orientation of $O$ yields no contradiction. 
The pencil is non-totally real if $2$ and $9$ lie together inside of the 
oval $2$ of $C_6$: a bad cubic is a cubic with an oval passing through $2$ and $9$, that is entirely contained in the empty oval $2$. To avoid this problem, let us choose the point $2$ {\em on\/} the corresponding empty oval $2$. This finishes the proof of ii). 
The last point iii) is an immediate consequence of ii).$\Box$


\vspace{3cm}

severine.fiedler@live.fr

\begin{thebibliography}{99}
\bibitem{F1}
T.~Fiedler: Eine Beschr\"ankung f\"ur die Lage von reellen algebraischen Kurven,
Beitr\"age Alg. Geom. No 11 (1981) 7-19
\bibitem{F2}
T.~Fiedler: Pencils of lines and the topology of real algebraic curves, Math. Izvestia
Vol. 2 (1983) No 1.
\bibitem{FLT1}
S.~Fiedler-Le Touz\'e: Orientations complexes des courbes alg\'ebriques
r\'eelles, These doctorale (2000).
\bibitem{FLT2}
S.~Fiedler-Le Touz\'e: Pencils of cubics with eight base points lying in convex position in $\mathbb{R}P^2$, arXiv[mathAG] 1012.2679, v2, 53 pages, Sept. 2012.
\bibitem{F-O}
S.~Fiedler-Le Touz\'e, S.~Orevkov: A flexible affine sextic which is
algebraically unrealizable, Journal of Algebraic Geometry {\bf 11},
 293-310 (2002).
\bibitem{G}
A.~Gabard: Ahlfors circle maps: historical ramblings, not yet on the web, available
at request asking the author alexandregabard@hotmail.com
\bibitem{K} 
V.M.~Kharlamov: Rigid isotopy clasification of real plane curves of degree 5, Funct.
Anal. Appl. 15 (1981)
\bibitem{K-V} 
V.M.~ Kharlamov, O.Ya.~ Viro: Extensions of the Gudkov-Rokhlin
congruence, Topology and Geometry Rokhlin seminar 1986, LNM 1346 Springer
357-406
\bibitem{N}
V.V.~Nikulin: Integer quadratic forms and some of their geometrical applications,
Math. USSR-Izv. 14 (1980), 103-167 
\bibitem{R} 
V. A.~Rokhlin: Complex topological characteristics of real algebraic curves, Russian Math. Surveys 33:5 (1978), 85-98 
\end{thebibliography}
\end{document}